\documentclass{article}

\usepackage[english]{babel}
\usepackage{amsmath, amssymb, amsthm, amsfonts, bm}
\usepackage{enumerate}
\usepackage{graphicx, color}
\usepackage[colorlinks=true]{hyperref}
\usepackage{enumitem}
\usepackage{pgf,tikz}
\usepackage{pgfplots}
\usepackage{cancel}

\numberwithin{equation}{section}
\newtheorem{theorem}{Theorem}[section]
\theoremstyle{definition}

\theoremstyle{remark}

\newcommand{\Rnum}[1]{\uppercase\expandafter{\romannumeral #1\relax}}

\newcommand{\mr}[1]{\mathrm{#1}}
\newcommand{\mb}[1]{\mathbb{#1}}
\newcommand{\mc}[1]{\mathcal{#1}}
\newcommand{\vsal}{\\\noalign{\vspace{0.2cm}}}

\DeclareMathOperator{\Av}{Av} % Average operator

\def\clap#1{\hbox to 0pt{\hss#1\hss}}

\title{On a constant related to the Bellman function of three integral variables of the dyadic maximal operator: Part B}
\author{Eleftherios N. Nikolidakis}
\date{\today}

\begin{document}
\maketitle

\begin{abstract}
We study the behaviour of the constant that is provided in the articles \cite{12} and \cite{13}, which is connected with the determination of the Bellman function of three integral variables of the dyadic maximal operator. More precisely we study the monotonicity properties of this constant with respect to the first variable from which it depends.
\end{abstract}

\section{Introduction} \label{sec:0}
The dyadic maximal operator on $\mb R^n$ is a useful tool in analysis and is defined by
\begin{equation} \label{eq:0p1}
	\mc M_d\varphi(x) = \sup\left\{ \frac{1}{|S|} \int_S |\varphi(u)|\,\mr du: x\in S,\ S\subseteq \mb R^n\ \text{is a dyadic cube} \right\},
\end{equation}
for every $\varphi\in L^1_\text{loc}(\mb R^n)$, where $|\cdot|$ denotes the Lebesgue measure on $\mb R^n$, and the dyadic cubes are those formed by the grids $2^{-N}\mb Z^n$, for $N=0, 1, 2, \ldots$.\\
It is well known that it satisfies the following weak type (1,1) inequality
\begin{equation} \label{eq:0p2}
	\left|\left\{ x\in\mb R^n: \mc M_d\varphi(x) > \lambda \right\}\right| \leq \frac{1}{\lambda} \int_{\left\{\mc M_d\varphi > \lambda\right\}} |\varphi(u)|\,\mr du,
\end{equation}
for every $\varphi\in L^1(\mb R^n)$, and every $\lambda>0$,
from which it is easy to get the following  $L^p$-inequality
\begin{equation} \label{eq:0p3}
	\|\mc M_d\varphi\|_p \leq \frac{p}{p-1} \|\varphi\|_p,
\end{equation}
for every $p>1$, and every $\varphi\in L^p(\mb R^n)$.
It is easy to see that the weak type inequality \eqref{eq:0p2} is the best possible. For refinements of this inequality one can consult \cite{6}.

It has also been proved that \eqref{eq:0p3} is best possible (see \cite{1} and \cite{2} for general martingales and \cite{21} for dyadic ones).
An approach for studying the behaviour of this maximal operator in more depth is the introduction of the so-called Bellman functions which play the role of generalized norms of $\mc M_d$. Such functions related to the $L^p$-inequality \eqref{eq:0p3} have been precisely identified in \cite{4}, \cite{5} and \cite{14}. For the study of the Bellman functions of $\mc M_d$, we use the notation $\Av_E(\psi)=\frac{1}{|E|} \int_E \psi$, whenever $E$ is a Lebesgue measurable subset of $\mb R^n$ of positive measure and $\psi$ is a real valued measurable function defined on $E$. We fix a dyadic cube  $Q$ and define the localized maximal operator $\mc M'_d\varphi$ as in \eqref{eq:0p1} but with the dyadic cubes $S$ being assumed to be contained in $Q$. Then for every $p>1$ we let
\begin{equation} \label{eq:0p4}
	B_p(f,F)=\sup\left\{ \frac{1}{|Q|} \int_Q (\mc M'_d\varphi)^p: \Av_Q(\varphi)=f,\ \Av_Q(\varphi^p)=F \right\},
\end{equation}
where $\varphi$ is nonnegative in $L^p(Q)$ and the variables $f, F$ satisfy $0<f^p\leq F$. By a scaling argument it is easy to see that \eqref{eq:0p4} is independent of the choice of $Q$ (so we may choose
$Q$ to be the unit cube $[0,1]^n$).
In \cite{5}, the function \eqref{eq:0p4} has been precisely identified for the first time. The proof has been given in a much more general setting of tree-like structures on probability spaces.

More precisely we consider a non-atomic probability space $(X,\mu)$ and let $\mc T$ be a family of measurable subsets of $X$, that has a tree-like structure similar to the one in the dyadic case (the exact definition can be seen in \cite{5}).
Then we define the dyadic maximal operator associated to $\mc T$, by
\begin{equation} \label{eq:0p5}
	\mc M_{\mc T}\varphi(x) = \sup \left\{ \frac{1}{\mu(I)} \int_I |\varphi|\,\mr \; d\mu: x\in I\in \mc T \right\},
\end{equation}
for every $\varphi\in L^1(X,\mu)$, $x\in X$.

This operator is related to the theory of martingales and satisfies essentially the same inequalities as $\mc M_d$ does. Now we define the corresponding Bellman function of four variables of $\mc M_{\mc T}$, by
\begin{multline} \label{eq:0p6}
	B_p^{\mc T}(f,F,L,k) = \sup \left\{ \int_K \left[ \max(\mc M_{\mc T}\varphi, L)\right]^p\mr \; d\mu: \varphi\geq 0, \int_X\varphi\,\mr \; d\mu=f, \right. \\  \left. \int_X\varphi^p\,\mr \; d\mu = F,\ K\subseteq X\ \text{measurable with}\ \mu(K)=k\right\},
\end{multline}
the variables $f, F, L, k$ satisfying $0<f^p\leq F $, $L\geq f$, $k\in (0,1]$.
The exact evaluation of \eqref{eq:0p6} is given in \cite{5}, for the cases where $k=1$ or $L=f$. In the first case, the author (in \cite{5}) precisely identifies the function $B_p^{\mc T}(f,F,L,1)$ by evaluating it in a first stage for the case where $L=f$. That is he precisely identifies $B_p^{\mc T}(f,F,f,1)$ (in fact $B_p^{\mc T}(f,F,f,1)=F \omega_p (\frac{f^p}{F})^p$, where                         $\omega_p: [0,1] \to [1,\frac{p}{p-1}]$ is the inverse function $H^{-1}_p$, of $H_p(z) = -(p-1)z^p + pz^{p-1}$). 

The proof of the above mentioned evaluation relies on a one-parameter integral inequality which is proved by arguments based on a linearization of the dyadic maximal operator. More precisely the author in \cite{5} proves that the inequality

\begin{equation}\label{eq:0p7}
	F\geq \frac{1}{(\beta+1)^{p-1}} f^p + \frac{(p-1)\beta}{(\beta+1)^p} \int_X (M_{\mathcal{T}}\varphi)^p \; d\mu,
\end{equation}
is true for every non-negative value of the parameter $\beta$ and sharp for one that depends on $f$, $F$ and $p$, namely for $\beta=\omega_p (\frac{f^p}{F})-1$. This gives as a consequence an upper bound for $B_p^{\mc T}(f,F,f,1)$, which after several technical considerations is proved to be best possible.Then  by using several calculus arguments the author in \cite{5} provides the evaluation of $B_p^{\mc T}(f,F,L,1)$ for every $L\geq f$. 

Now in \cite{14} the authors give a direct proof of the evaluation of $B_p^{\mc T}(f,F,L,1)$ by using alternative methods. Moreover in the second case, where $L=f$, the author (in \cite{5}) uses the evaluation of $B_p^{\mc T}(f,F,f,1)$ and provides the evaluation of the more general $B_p^{\mc T}(f,F,f,k)$, $k\in (0,1]$.

Our aim in this article is to study further the results of \cite{12} and \cite{13} in order to approach the following Bellman function problem (of three integral variables)

\begin{multline} \label{eq:0p8}
	B_{p,q}^{\mc T}(f,A,F) = \sup \left\{ \int_X \left(\mc M_{\mc T}\varphi\right)^p\mr \; d\mu: \varphi\geq 0, \int_X\varphi\,\mr \; d\mu=f, \right. \\  \left. \int_X\varphi^q\,\mr \; d\mu = A,\ \int_X\varphi^p\,\mr \; d\mu = F\right\},
\end{multline}
where $1<q<p$, and the variables $f,A,F$ lie in the domain of definition of the above problem. Certain progress for the above problem can be seen in \cite{11}.

In \cite{12} we have proved that whenever $0<\frac{x^q}{\kappa^{q-1}}<y\leq x^{\frac{p-q}{p-1}}\cdot z^{\frac{q-1}{p-1}}\;\Leftrightarrow\; 0<s_1^{\frac{q-1}{p-1}}\leq s_2<1$, (where $s_1,s_2$ are defined below), there is a constant $t=t(s_1,s_2)$ for which if $h:(0,\kappa]\longrightarrow\mb{R}^{+}$ satisfies  $\int_{0}^{\kappa}h=x$\,,\; $\int_{0}^{\kappa}h^q=y$ and $\int_{0}^{\kappa}h^p=z$ then

	$$\int_{0}^{\kappa}\bigg(\frac{1}{t}\int_{0}^{t}h\bigg)^p dt\leq t^p(s_1,s_2)\cdot\int_{0}^{\kappa}h^p.$$

More precisely $t(s_1,s_2)=t$ is the greatest element of $\big[{1,t(0)}\big]$ for which $F_{s_1,s_2}(t)\leq 0$ where $F_{s_1,s_2}$ is defined in \cite{12}. Moreover for each such fixed $s_1,s_2$ 
\[t=t(s_1,s_2)=\min\Big\{{t(\beta)\;:\;\beta\in\big[0,\tfrac{1}{p-1}\big]}\Big\}\]
where $t(\beta)=t(\beta,s_1,s_2)$ is defined in \cite{12}. That is we find a constant $t=t(s_1,s_2)$ for which the above inequality  is satisfied for all $h:(0,\kappa]\longrightarrow\mb{R}^{+}$ as mentioned above. Note that $s_1,s_2$ depend by a certain way on $x,y,z$, namely $s_1=\frac{x^p}{\kappa^{p-1}z}$,  $s_2=\frac{x^q}{\kappa^{q-1}y}$ and $F_{s_1,s_2}(\cdot)$ is given in terms of $s_1,s_2$.

In this article we study the monotonicity behaviour of the function 
$t(s_1,s_2)$ with respect to the first variable $s_1$. This may enable us in the future to determine \eqref{eq:0p6}, by using also the results in \cite{11}, \cite{12} and \cite{13}.

We need to mention that the extremizers for the standard Bellman function $B_p^{\mc T}(f,F,f,1)$ have been studied in \cite{7}, and in \cite{9} for the case $0<p<1$. Also in \cite{8} the extremal sequences of functions for the respective Hardy operator problem have been studied.  Additionally further study of the dyadic maximal operator can be seen in \cite{10,14} where symmetrization principles for this operator are presented, while other approaches for the determination of certain Bellman functions are given in \cite{16,17,18,19,20}. Moreover results related to applications of dyadic maximal operators can be seen in \cite{15}.

\bigskip

\section{Monotonicity of $t(s_1,s_2)$ with respect to the variable $s_1$}\label{sec:3}

As we have already proved whenever $t(s_1,s_2)$ is given by $t_{s_1,s_2}(0)$, that is whenever $(s_1,s_2)$ is such that $\min\big\{{t(s_1,s_2,\gamma)\;:\; \gamma\in\big[{0,\frac{1}{p-1}}\big]}\big\}=t(s_1,s_2,0)=t_{s_1,s_2}(0)$, the following is true $\frac{\partial t_{s_1,s_2}(0)}{\partial s_1}=\frac{\partial t(0)}{\partial s_1}<0$.

\medskip

We now investigate the rest domain of definition of $(s_1,s_2)$. Thus we assume that the variables $(s_1,s_2)$ satisfy 
\begin{align}
	0<s_1<\min\Big\{{s_2^{\frac{p-1}{q-1}},s_1^{(0)}=\big({\tfrac{1}{h({s_2})}}\big)^{\frac{1}{q}}}\Big\}\,.\label{eq:3p1}
\end{align}
For any such $(s_1,s_2)$, it is previously seen that $t({s_1,s_2})=t$ is given by the following equation:
\begin{align}
	q\Big({p\,\omega_q(\tau)^{q-1}-({p-1})\,\omega_q(\tau)^{q}}\Big)\bigg({t^{p-q}-\frac{s_1}{s_2}}\bigg)=(p-q)s_1\cdot\alpha(s_2)\,,\label{eq:3p2}
\end{align}
where $\alpha(s_2)=\frac{\omega_q(s_2)^q}{s_2}-1$\,,\; $s_2\in(0,1)$, and $\tau(s_1,s_2,t)=\tau=\frac{p-q}{p}\,\frac{t^p-s_1}{t^{p-q}-\frac{s_1}{s_2}}$\,.

We will prove that for $(s_1,s_2)$ as in \eqref{eq:3p1}, the function $t(s_1,s_2)$ is strictly decreasing with respect to the variable $s_1$. We first evaluate $\frac{\partial \tau}{\partial s_1}$ :
\begin{align}
	\frac{\partial \tau}{\partial s_1}&=\frac{p-q}{p}\,\frac{\big({p\,\frac{\partial t}{\partial s_1}\,t^{p-1}-1}\big)\big({t^{p-q}-\frac{s_1}{s_2}}\big)-(t^p-s_1)\big({(p-q)\,\frac{\partial t}{\partial s_1}\,t^{p-q-1}-\frac{1}{s_2}}\big)}{\big({t^{p-q}-\frac{s_1}{s_2}}\big)^2}\nonumber\vsal
	&=\frac{p-q}{p}\,\frac{1}{\big({t^{p-q}-\frac{s_1}{s_2}}\big)^2}\bigg[p\,t^{2p-q-1}\frac{\partial t}{\partial s_1} -p\,\frac{\partial t}{\partial s_1}\,t^{p-1}\frac{s_1}{s_2}-t^{p-q}+\frac{s_1}{s_2}\,-\nonumber\vsal  &\hspace{2.0cm}- (p-q)\,t^{2p-q-1}\frac{\partial t}{\partial s_1}+\frac{t^p}{s_2}+(p-q)\,\frac{\partial t}{\partial s_1}\,s_1\,t^{p-q-1}-\frac{s_1}{s_2}\bigg]\nonumber\vsal
	&=\frac{p-q}{p}\,\frac{1}{\big({t^{p-q}-\frac{s_1}{s_2}}\big)^2}\bigg[q\,t^{2p-q-1}\frac{\partial t}{\partial s_1} -p\,\frac{\partial t}{\partial s_1}\,t^{p-1}\frac{s_1}{s_2}-t^{p-q}+\frac{t^p}{s_2}\,+\nonumber\vsal
	&\hspace{6.9cm} +(p-q)\,\frac{\partial t}{\partial s_1}\,s_1\,t^{p-q-1}\bigg]\;\;\Rightarrow\nonumber
\end{align}

\begin{align}
	\frac{\partial \tau}{\partial s_1}&=\frac{p-q}{p}\,\frac{1}{\big({t^{p-q}-\frac{s_1}{s_2}}\big)^2}\Bigg[\frac{\partial t}{\partial s_1}\,t^{p-q-1}\bigg({q\,t^p-p\,t^q\,\frac{s_1}{s_2}+(p-q)s_1}\bigg)\,+\nonumber\vsal
	&\hspace{7.cm} +\bigg({\frac{t^p}{s_2}-t^{p-q}}\bigg)\Bigg]\,.\label{eq:3p3}
\end{align}
We now differentiate \eqref{eq:3p2} with respect to $s_1$ and we get:

\begin{align*}
	\frac{p-q}{q}\,\alpha(s_2)&=\bigg({(p-q)\,t^{p-q-1}\frac{\partial t}{\partial s_1}-\frac{1}{s_2}}\bigg)\Big({p\,\omega_q(\tau)^{q-1}-({p-1})\,\omega_q(\tau)^{q}}\Big)\,+\vsal
	&\hspace{1.5cm} +\bigg({t^{p-q}-\frac{s_1}{s_2}}\bigg)\bigg[p(q-1)\,\omega_q(\tau)^{q-2}\omega'_q(\tau)\,\frac{\partial \tau}{\partial s_1}\,-\vsal
	&\hspace{4.5cm}-({p-1})q\,\omega_q(\tau)^{q-1}\omega'_q(\tau)\,\frac{\partial \tau}{\partial s_1}\bigg]\vsal
	&=\bigg({(p-q)\,t^{p-q-1}\frac{\partial t}{\partial s_1}-\frac{1}{s_2}}\bigg)\Big({p\,\omega_q(\tau)^{q-1}-({p-1})\,\omega_q(\tau)^{q}}\Big)\,+\vsal
	&\qquad+\bigg({t^{p-q}-\frac{s_1}{s_2}}\bigg)\,\omega'_q(\tau)\,\frac{\partial \tau}{\partial s_1}\,\omega_q(\tau)^{q-2}\big({p(q-1)-(p-1)q\,\omega_q(\tau)}\big)\vsal
	&\overset{\eqref{eq:3p2}}{\underset{\eqref{eq:3p3}}{=\!=\!=}}\bigg({(p-q)\,t^{p-q-1}\frac{\partial t}{\partial s_1}-\frac{1}{s_2}}\bigg)\,\frac{p-q}{q}\,\frac{s_1\,\alpha(s_2)}{t^{p-q}-\frac{s_1}{s_2}}\,+\vsal
	&\qquad +\bigg({t^{p-q}-\frac{s_1}{s_2}}\bigg)\,\frac{1}{q(q-1)\,\omega_q(\tau)^{q-2}\big({1-\omega_q(\tau)}\big)}\,\omega_q(\tau)^{q-2}\,\cdot\vsal
	&\hspace{2.cm}\cdot\big({p(q-1)-(p-1)q\,\omega_q(\tau)}\big)\,\frac{p-q}{p}\,\frac{1}{\big({t^{p-q}-\frac{s_1}{s_2}}\big)^2}\,\cdot\vsal
	&\hspace{1.15cm}\cdot \bigg[{\frac{\partial t}{\partial s_1}\,t^{p-q-1}\bigg({q\,t^p-p\,t^q\frac{s_1}{s_2}+(p-q)s_1}\bigg)+\frac{t^p}{s_2}-t^{p-q}}\bigg]\;\;\Rightarrow\vsal
	\alpha(s_2)&=\bigg({(p-q)\,t^{p-q-1}\frac{\partial t}{\partial s_1}-\frac{1}{s_2}}\bigg)\,\frac{s_1\,\alpha(s_2)}{t^{p-q}-\frac{s_1}{s_2}} \,+\vsal
	&\qquad+\frac{q}{p-q}\bigg({t^{p-q}-\frac{s_1}{s_2}}\bigg)\,\frac{1}{q(q-1)\,\omega_q(\tau)^{q-2}\big({1-\omega_q(\tau)}\big)}\,\omega_q(\tau)^{q-2}\,\cdot
\end{align*}

\pagebreak

\begin{align*}
	&\hspace{1.2cm}\cdot\big({p(q-1)-(p-1)q\,\omega_q(\tau)}\big)\,\frac{p-q}{p}\,\frac{1}{\big({t^{p-q}-\frac{s_1}{s_2}}\big)^2}\,\cdot\vsal &\hspace{2.5cm}\cdot\bigg[{\frac{\partial t}{\partial s_1}\,t^{p-q-1}
		\bigg({q\,t^p-p\,t^q\frac{s_1}{s_2}+(p-q)s_1}\bigg)+\frac{t^p}{s_2}-t^{p-q}}\bigg]	\vsal
	&=\frac{s_1\,\alpha(s_2)}{t^{p-q}-\frac{s_1}{s_2}}\bigg({(p-q)\,t^{p-q-1}\frac{\partial t}{\partial s_1}-\frac{1}{s_2}}\bigg)\,+\vsal
	&\hspace{1.2cm}+\frac{1}{(q-1)p}\,\frac{1}{\omega_q(\tau)-1}\cdot\big({(p-1)q\,\omega_q(\tau)-p(q-1)}\big)\,\frac{1}{t^{p-q}-\frac{s_1}{s_2}}\,\cdot\vsal
	&\hspace{2.4cm}\cdot\bigg[\frac{\partial t}{\partial s_1}\,t^{p-q-1}
	\bigg({q\,t^p-p\,t^q\frac{s_1}{s_2}+(p-q)s_1}\bigg)+\frac{t^p}{s_2}-t^{p-q} \bigg]\;\;\Rightarrow
\end{align*}
\begin{align*}
	\bigg({t^{p-q}-\frac{s_1}{s_2}}\bigg)\,\alpha(s_2)&=s_1\,\alpha(s_2)\bigg({(p-q)\,t^{p-q-1}\frac{\partial t}{\partial s_1}-\frac{1}{s_2}}\bigg)\,+\vsal
	&\qquad+\frac{1}{p(q-1)}\,\frac{(p-1)q\,\omega_q(\tau)-p(q-1)}{\omega_q(\tau)-1}\,t^{p-q-1}\,\frac{\partial t}{\partial s_1}\,\cdot\vsal
	&\hspace{4.cm}\cdot \bigg({q\,t^p-p\,t^q\frac{s_1}{s_2}+(p-q)s_1}\bigg)\,+\vsal
	&\hspace{0.8cm}+\frac{1}{p(q-1)}\,\frac{(p-1)q\,\omega_q(\tau)-p(q-1)}{\omega_q(\tau)-1}\bigg({\frac{t^p}{s_2}-t^{p-q} }\bigg)\quad\Rightarrow
\end{align*}
\begin{align*}
	&\bigg\{ \frac{1}{p(q-1)}\,\frac{(p-1)q\,\omega_q(\tau)-p(q-1)}{\omega_q(\tau)-1}\,t^{p-q-1}\bigg({q\,t^p-p\,t^q\frac{s_1}{s_2}+(p-q)s_1}\bigg)\,+\nonumber\vsal
	&\hspace{6.75cm}+({p-q})s_1\,\alpha(s_2)\,t^{p-q-1}\bigg\}\,\frac{\partial t}{\partial s_1}\nonumber\vsal
	&=\bigg({t^{p-q}-\frac{s_1}{s_2}}\bigg)\,\alpha(s_2)+\frac{s_1}{s_2}\,\alpha(s_2)\,-\nonumber\vsal
	&\hspace{2.75cm}-\frac{1}{p(q-1)}\,\frac{(p-1)q\,\omega_q(\tau)-p(q-1)}{\omega_q(\tau)-1}\,\bigg({\frac{t^p}{s_2}-t^{p-q}}\bigg)\quad\Rightarrow\vsal
	&t^{p-q-1}\frac{\partial t}{\partial s_1}\bigg\{\frac{1}{p(q-1)}\,\frac{(p-1)q\,\omega_q(\tau)-p(q-1)}{\omega_q(\tau)-1}\bigg({q\,t^{p}-p\,t^q\,\frac{s_1}{s_2}+(p-q)s_1}\bigg)\,+ \nonumber\vsal
	&\hspace{8.5cm} + (p-q)s_1\,\alpha(s_2) \bigg\}\nonumber
\end{align*}
\begin{align}
	&=t^{p-q}\,\alpha(s_2)-\frac{1}{p(q-1)}\,\frac{(p-1)q\,\omega_q(\tau)-p(q-1)}{\omega_q(\tau)-1}\,t^{p-q}\bigg({\frac{t^q}{s_2}-1}\bigg)\hspace{0.85cm}\Rightarrow\nonumber\vsal
	&\hspace{4.0cm}\frac{\partial t}{\partial s_1}\cdot\delta(s_1,s_2)=t\,\gamma(s_1,s_2)\,,\label{eq:3p4}
\end{align}
where 
\begin{align}
	\delta(s_1,s_2)&:=\frac{1}{p(q-1)}\,\frac{(p-1)q\,\omega_q(\tau)-p(q-1)}{\omega_q(\tau)-1}\bigg({q\,t^{p}-p\,t^q\,\frac{s_1}{s_2}+(p-q)s_1}\bigg)
	\,+\nonumber\vsal
	&\hspace{6.5cm}+ (p-q)s_1\,\alpha(s_2)\label{eq:3p5}
\end{align}
and 
\begin{align}
	\gamma(s_1,s_2)&:=\alpha(s_2)-\frac{1}{p(q-1)}\,\frac{(p-1)q\,\omega_q(\tau)-p(q-1)}{\omega_q(\tau)-1}\bigg({\frac{t^q}{s_2}-1}\bigg)\,.\label{eq:3p6}
\end{align}

We now study the quantity 
\[\lambda(t)=q\,t^{p}-p\,t^q\,\frac{s_1}{s_2}+(p-q)s_1\,.\]
We prove that $\lambda(t)\geq \lambda(1)>0$, \; $\forall\,t\in \big({1,\frac{p}{p-1}}\big)$. We have $\lambda'(t)=p\,q\,t^{p-1}-p\,q\,t^{q-1}\,\frac{s_1}{s_2}\cong t^{p-q}-\frac{s_1}{s_2}>1-1=0$. Thus for every $t\in \big[{1,\frac{p}{p-1}}\big)$ we have $\lambda(t)\geq \lambda(1)$.

We now prove the second assertion that is made above, namely that:
\begin{align*}
	q-p\,\frac{s_1}{s_2}+(p-q)s_1=\lambda(1)>0\quad&\Leftrightarrow\vsal
	(p-q)s_1>p\,\frac{s_1}{s_2}-q\quad&\Leftrightarrow\quad p\,\frac{s_1}{s_2}<(p-q)s_1+q\,.
\end{align*}
For this last inequality it is enough to prove that $p\,\frac{s_1}{s_1^{\frac{q-1}{p-1}}}<(p-q)s_1+q$, since $s_2\geq s_1^{\frac{q-1}{p-1}}$ in our domain of definition. Thus it is enough to prove 
\begin{align*}
	p\,s_1^{\frac{p-q}{p-1}}<(p-q)s_1+q\,,\;\;  \forall\,s_1\in({0,1})\,.\qquad(*)
\end{align*}
We set $\varphi(s_1):=p\,s_1^{\frac{p-q}{p-1}}-(p-q)s_1-q$. Then we notice that $\varphi(1)=0$, while 
\begin{align*}
	\varphi'(s_1)=p\,\frac{p-q}{p-1}\,s_1^{-\frac{q-1}{p-1}}-(p-q)\cong \frac{p}{p-1}\,s_1^{-\frac{q-1}{p-1}}-1\stackrel{s_1<1}{>}0\,.
\end{align*}
Thus $\varphi(s_1)\stackrel{s_1<1}{<}\varphi(1)=0$, whenever $(*)$ is true. Thus from our conclusion we get $\lambda(t)\geq \lambda(1)>0$, \; $\forall\,t\in \big[{1,\frac{p}{p-1}}\big)$, where 
\[\lambda(t)=q\,t^{p}-p\,t^q\,\frac{s_1}{s_2}+(p-q)s_1\,.\]

We now continue our study for the sign of  $\frac{\partial t}{\partial s_1}$:
For the moment, fix $(s_1,s_2)$ such that: $\omega_q(s_2)=\omega_p(s_1)\;\Leftrightarrow\; s_2=H_q({\omega_p(s_1)})$. Note that such a point $(s_1,s_2)$ belongs to the domain of definition which is described by \eqref{eq:3p1}. We prove that for such a point we have: $\gamma(s_1,s_2)<0$ where $\gamma(s_1,s_2)$ is given by \eqref{eq:3p6}. 
That is 
\begin{align*}
	\gamma(s_1,s_2)&=\alpha(s_2)-\frac{1}{p(q-1)}\,\frac{(p-1)q\,\omega_q(\tau)-p(q-1)}{\omega_q(\tau)-1}\bigg({\frac{t^q}{s_2}-1}\bigg)\,.
\end{align*}
We set $\omega_p(s_1)=\omega_q(s_2)=\lambda$. Then as we have seen before $t(s_1,s_2)=\omega_p(s_1)=\lambda$.  Thus we have $\tau:=\frac{p-q}{p}\,\frac{t^p-s_1}{t^{p-q}-\frac{s_1}{s_2}}=s_2$. Indeed, it's enough to prove that:
\begin{align*}
	\frac{p-q}{p}\,\frac{\omega_p(s_1)^p-s_1}{\omega_p(s_1)^{p-q}-\frac{s_1}{H_q(\omega_p(s_1))}}&=s_2\;\big({=H_q(\omega_p(s_1))}\big)&&\Leftrightarrow\vsal
	\frac{p-q}{p}\,\frac{\lambda^p-H_p(\lambda)}{\lambda^{p-q}-\frac{H_p(\lambda)}{H_q(\lambda)}}&=H_q(\lambda)&&\Leftrightarrow\vsal
	\frac{p-q}{p}\,\big({\lambda^p-H_p(\lambda)}\big)&=\lambda^{p-q}H_q(\lambda)-H_p(\lambda)&&\Leftrightarrow\vsal
	\frac{p-q}{p}\,\big({\lambda^p-p\lambda^{p-1}+(p-1)\lambda^p}\big)&=\lambda^{p-q}\big({q\lambda^{q-1}-(q-1)\lambda^{q}}\big)\,-\vsal
	&\hspace{1.3cm}-\big({p\lambda^{p-1}-(p-1)\lambda^{p}}\big)&&\Leftrightarrow\vsal
	(p-q)\big({\lambda^p-\lambda^{p-1}}\big)&=\lambda^{p-1}({q-(q-1)\lambda})-\lambda^{p-1}({p-(p-1)\lambda})&&\Leftrightarrow\vsal
	(p-q)(\lambda-1)&=q-(q-1)\lambda-p+(p-1)\lambda&&\Leftrightarrow\vsal
	(p-q)(\lambda-1)&=(p-q)\lambda-(p-q)\,,
\end{align*}
which is obviously true. We have just proved that whenever $\omega_p(s_1)=\omega_q(s_2)$, then $\tau=s_2$. 

But then 
\begin{align*}
	\gamma(s_1,s_2)&=\bigg({\frac{\omega_q(s_2)^q}{s_2}-1}\bigg)-\frac{1}{p(q-1)}\,\frac{(p-1)q\,\omega_q(s_2)-p(q-1)}{\omega_q(s_2)-1}\bigg({\frac{\omega_q(s_2)^q}{s_2}-1}\bigg)\vsal
	&\stackrel{\omega_q(s_2)=\lambda}{=\!=\!=\!=\!=\!=}\bigg({\frac{\lambda^q}{H_q(\lambda)}-1}\bigg)\bigg({1-\frac{1}{p(q-1)}\,\frac{(p-1)q\lambda-p(q-1)}{\lambda-1}}\bigg)\vsal
	&\stackrel{\text{equall  sign}}{\cong}p(q-1)({\lambda-1})-(p-1)q\lambda+p(q-1)\vsal
	&=p(q-1)\lambda-(p-1)q\lambda\vsal
	&=(q-p)\lambda=-(p-q)\lambda<0\,.
\end{align*}
Moreover since $q\,t^p-p\,t^q\,\frac{s_1}{s_2}+(p-q)s_1$ is always strictly positive, we get by \eqref{eq:3p5} that $\delta(s_1,s_2)>0$, for every $(s_1,s_2)$. Thus we have proved that $\gamma(s_1,s_2)<0$, whenever $ s_2=H_q({\omega_p(s_1)})$.

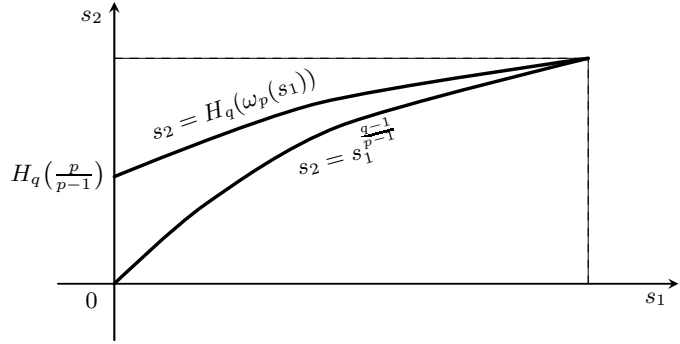
\begin{center}
	\begin{tikzpicture}[line cap=round,line join=round,>=stealth,x=2.0cm,y=1.0cm,scale=1.5]
		\draw[->,color=black,line width=0.7pt] (-0.25,0.) -- (2.5,0.);
		\foreach \x in {}
		\draw[shift={(\x,0)},color=black] (0pt,-2.5pt) -- (0pt,2.5pt)
		node[below] at (0,-1pt) {};	
		%\draw[color=black] (4,0) node[above]  {{\small{$y$}}};
		\draw[->,color=black,line width=0.8pt] (0,-0.5) -- (0.,2.5);
		\foreach \y in {}
		\draw[shift={(0,\y)},color=black] (-2pt,0pt) -- (2pt,0pt) node[left] at (-1pt,0pt) {};
		\draw [fill=black,dashed] (0,2) -- (2.1,2);
		\draw [fill=black,dashed] (2.1,0) -- (2.1,2);
		\begin{small}
			\draw[color=black] (0.5,1.75) node[below,rotate=19] {$s_2=H_q({\omega_p(s_1)})$};
			\draw[color=black] (0,0.95) node[left] {$H_q\big({\frac{p}{p-1}}\big)$};
			\draw[color=black] (1.,1.45) node[below,rotate=15] {$s_2=s_1^{\frac{q-1}{p-1}}$};
			\draw[color=black] (-0.1,0) node[below] {$0$};
			\draw[color=black] (2.4,0) node[below]  {$s_1$};
			\draw[color=black] (-0.1,2.5) node[below]  {$s_2$};
\draw [line width=1.3pt ] plot [smooth,samples=200] coordinates {(0.,0.95) (0.9,1.6)(2.1,2)};
\draw [line width=1.3pt] plot [smooth,samples=200] coordinates {(0.,0) (0.4, 0.7) (1,1.4) (2.1,2)};
		\end{small}
	\end{tikzpicture}
\end{center}

Fix now $s_2<H_q\big({\frac{p}{p-1}}\big)=H_q({\omega_p(0)})$. At this point note that the facts mentioned below hold by using standard arguments involving limits of sequences and subsequences. Since $t({s_1,s_2})=t$ is bounded ($1\leq t\leq \frac{p}{p-1}$), by letting $s_1\to 0^{+}$ (with $s_2<H_q\big({\frac{p}{p-1}}\big)$ being fixed) we obtain that $t({s_1,s_2})\to t_0$ for some $t_0\in\big[{1,\frac{p}{p-1}}\big]$. Thus by \eqref{eq:3p1} we get that 
\begin{align*}
&p\,\omega_q(\tau)^{q-1}-(p-1)\,\omega_q(\tau)^{q}\xrightarrow[s_2\,\text{fixed},\, s_2<H_q({\frac{p}{p-1}})]{s_1\,\cong\, 0^{+}}0\quad\Rightarrow\vsal
&\omega_q(\tau)\cong \frac{p}{p-1}\quad\Rightarrow\quad \tau\cong H_q\big({\tfrac{p}{p-1}}\big)\,.
\end{align*}
But also $\tau=\frac{p-q}{p}\,\frac{t^p-s_1}{t^{p-q}-\frac{s_1}{s_2}}\cong \frac{p-q}{p}\,t_0^q$. Thus we must have that: 
\begin{align*}
\frac{p-q}{p}\,t_0^q&=H_q\big({\tfrac{p}{p-1}}\big)&&\Rightarrow\vsal
t_0^q&=\bigg({q\Big({\frac{p}{p-1}}\Big)^{q-1}-(q-1)\Big({\frac{p}{p-1}}\Big)^q}\bigg)\frac{p}{p-q}\vsal
  &=\Big({\frac{p}{p-1}}\Big)^{q-1}\bigg({q-(q-1)\frac{p}{p-1}}\bigg)\,\frac{p}{p-q}\vsal
  &=\Big({\frac{p}{p-1}}\Big)^{q-1}\,\frac{p}{p-q}\,\frac{qp-q-qp+p}{p-1}=\Big({\frac{p}{p-1}}\Big)^{q} &&\Rightarrow\vsal
t_0&=\frac{p}{p-1}\,.
\end{align*}
Also $\tau\cong H_q\big({\tfrac{p}{p-1}}\big)\;\Rightarrow\; \omega_q(\tau)\cong \tfrac{p}{p-1}$\,. Thus for any fixed $s_2<H_q\big({\tfrac{p}{p-1}}\big)$, if we compute $\lim_{s_1\to 0^{+}}\gamma(s_1,s_2):=\gamma_0$, we have 
\begin{align*}
\gamma_0&=\alpha(s_2)-\frac{1}{p(q-1)}\,\frac{(p-1)q\,\frac{p}{p-1}-p(q-1)}{\frac{p}{p-1}-1}\Big({\frac{t_0^q}{s_2}-1}\Big)\vsal
&=\alpha(s_2)-\frac{1}{p(q-1)}\,\frac{p}{\frac{1}{p-1}}\Bigg({\frac{\big({\tfrac{p}{p-1}}\big)^q}{s_2}-1}\Bigg)&&\Rightarrow
\end{align*}

\begin{align*}
\gamma_0&=\alpha(s_2)-\frac{p-1}{q-1}\Bigg({\frac{\big({\tfrac{p}{p-1}}\big)^q}{s_2}-1}\Bigg)\vsal
        &=\bigg({\frac{\omega_q(s_2)^q}{s_2}-1}\bigg)-\frac{p-1}{q-1}\Big({\frac{p}{p-1}}\Big)^q\frac{1}{s_2}+\frac{p-1}{q-1}\vsal
      &=\frac{\omega_q(s_2)^q-\alpha}{s_2}-1+\frac{p-1}{q-1}:=F(s_2)\,,
 \end{align*}
where $\alpha=\frac{p-1}{q-1}\big({\tfrac{p}{p-1}}\big)^q$. Then $F(s_2)$ defined right above satisfies 
\begin{align}
F({s_2})<\frac{\big({\tfrac{q}{q-1}}\big)^q-\frac{p-1}{q-1}\big({\tfrac{p}{p-1}}\big)^q}{s_2}-1+\frac{p-1}{q-1}\,.\label{eq:3p7}
\end{align}
We now prove the following inequalities:
\begin{align*}
\Big({\frac{q}{q-1}}\Big)^q<\frac{p-1}{q-1}\Big({\frac{p}{p-1}}\Big)^q
=:\alpha\quad\Leftrightarrow\quad \frac{p^q}{(p-1)^{q-1}}>\frac{q^q}{(q-1)^{q-1}}\,.
\end{align*}
For this purpose we define the function $\varphi_q(x)=\frac{x^q}{({x-1})^{q-1}}$,\;$x>1$, whenever $q>1$. It is enough then to prove $\varphi'_q(x)>0$, $\forall\,x>q$. Indeed $\varphi'_q(x)\cong q\,x^{q-1}(x-1)^{q-1}-x^q(q-1)(x-1)^{q-2}\cong q(x-1)-(q-1)\,x=x-q>0$, $\forall\,x>q$. Thus since 
\begin{align*}
p>q>1\;\Rightarrow\;\frac{p^q}{(p-1)^{q-1}}>\frac{q^q}{(q-1)^{q-1}}\;\Rightarrow\; \Big({\frac{q}{q-1}}\Big)^q-\frac{p-1}{q-1}\Big({\frac{p}{p-1}}\Big)^q<0\,,
\end{align*}
then by \eqref{eq:3p7} $\lim_{s_2\to 0^{+}}F(s_2)=-\infty$. 

We wish to prove that $F(s_2)<0$, for any $0<s_2<H_q\big({\tfrac{p}{p-1}}\big)$. Remember that $F(s_2):=\frac{\omega_q(s_2)^q-\alpha}{s_2}-1+\frac{p-1}{q-1}$, where $\alpha=\frac{p-1}{q-1}\big({\tfrac{p}{p-1}}\big)^q$. Note that:

\begin{align}
F'(s_2)&=\frac{q\,\omega_q(s_2)^{q-1}\,\frac{1}{q(q-1)\omega_q(s_2)^{q-2}(1-\omega_q(s_2))}s_2-\omega_q(s_2)^q}{s_2^2}+\frac{\alpha}{s_2^2}\nonumber\vsal
  &\cong-\underbrace{\frac{\omega_q(s_2)}{({q-1})\big({\omega_q(s_2)-1}\big)}\,s_2-\omega_q(s_2)^q}\limits_{=\,-G(s_2)}+\alpha=-G(s_2)+\alpha\,.\label{eq:3p8}
\end{align}
As we shall see latter the function $G(s_2)$,\;$s_2\in(0,1)$ defined right above, is \emph{strictly increasing}, thus since $0<s_2<H_q\big({\tfrac{p}{p-1}}\big)$ we have 
\begin{align*}
G(s_2)<G\big({H_q\big({\tfrac{p}{p-1}}\big)}\big)=\frac{1}{q-1}\,\frac{\frac{p}{p-1}}{\frac{1}{p-1}}\,H_q\big({\tfrac{p}{p-1}}\big)+\Big({\frac{p}{p-1}}\Big)^q\,.
\end{align*}
Thus by \eqref{eq:3p8} above we get:
\begin{align*}
F'({s_2})& > -\bigg[{\frac{p}{q-1}\,H_q\big({\tfrac{p}{p-1}}\big)+\Big({\frac{p}{p-1}}\Big)^q}\bigg]+\alpha\vsal
 &=-\frac{p}{q-1}\bigg({\frac{p-q}{p}\Big({\frac{p}{p-1}}\Big)^q}\bigg)+\Big({\frac{p}{p-1}}\Big)^q+\frac{p-1}{q-1}\Big({\frac{p}{p-1}}\Big)^q\vsal
 &=\Big({\frac{p}{p-1}}\Big)^q+\frac{p-1}{q-1}\Big({\frac{p}{p-1}}\Big)^q-\frac{p-q}{q-1}\Big({\frac{p}{p-1}}\Big)^q\vsal
 &\cong  1+\frac{p-1}{q-1}-\frac{p-q}{q-1}=0\,.
\end{align*}
Thus $F'({s_2}) > 0$, \;$\forall\,0<s_2<H_q\big({\tfrac{p}{p-1}}\big)$\,. But also 
\begin{align*}
F\big({H_q\big({\tfrac{p}{p-1}}\big)}\big)&=\frac{\omega_q\big({H_q\big({\tfrac{p}{p-1}}\big)}\big)^q}{H_q\big({\tfrac{p}{p-1}}\big)}-\frac{\frac{p-1}{q-1}\,\big({\tfrac{p}{p-1}}\big)^q}{H_q\big({\tfrac{p}{p-1}}\big)}-1+\frac{p-1}{q-1}\vsal
&=\frac{\big({\tfrac{p}{p-1}}\big)^q}{\frac{p-q}{p}\big({\tfrac{p}{p-1}}\big)^q}-\frac{\frac{p-1}{q-1}\,\big({\tfrac{p}{p-1}}\big)^q}{\frac{p-q}{p}\big({\tfrac{p}{p-1}}\big)^q}-1+\frac{p-1}{q-1}\vsal
&=\frac{p}{p-q}-\frac{p}{p-q}\,\frac{p-1}{q-1}-1+\frac{p-1}{q-1}\vsal
&=\frac{p}{p-q}\bigg({1-\frac{p-1}{q-1}}\bigg)-\bigg({1-\frac{p-1}{q-1}}\bigg)\vsal
&=\bigg({\frac{p}{p-q}-1}\bigg)\bigg({1-\frac{p-1}{q-1}}\bigg)<0\,.
\end{align*}
Thus $F(s_2)<0$,\;$\forall\,s_2\in\big({0,H_q\big({\tfrac{p}{p-1}}\big)}\big)$.

We have just proved that whenever $s_2$ satisfies $0<s_2<H_q\big({\tfrac{p}{p-1}}\big)$ then $\lim_{s_1\to 0^{+}}\gamma(s_1,s_2)<0$. 

We continue proving that whenever $s_2$ satisfies $H_q\big({\tfrac{p}{p-1}}\big)\leq s_2<1$ then also have $\gamma(s_1,s_2)\underset{s_1\to 0^{+}}{\cong}$ strictly negative number.

Remember that 
\begin{align*}
	\gamma(s_1,s_2)&=\alpha(s_2)-\frac{1}{p(q-1)}\,\frac{(p-1)q\,\omega_q(\tau)-p(q-1)}{\omega_q(\tau)-1}\bigg({\frac{t^q}{s_2}-1}\bigg)\,.
\end{align*}
 by \eqref{eq:3p6}. As before (in the previous case) we can see that  again we have $t(s_1,s_2)\underset{s_1\to 0^{+}}{\cong}\frac{p}{p-1}$ and $\tau\underset{s_1\to 0^{+}}{\cong}H_q\big({\tfrac{p}{p-1}}\big)$. Thus: 

\begin{align}
\gamma(s_1,s_2)&\cong\alpha(s_2)-\frac{1}{p(q-1)}\,\frac{(p-1)q\,\frac{p}{p-1}-p(q-1)}{\frac{p}{p-1}-1}\Bigg({\frac{\big({\tfrac{p}{p-1}}\big)^q}{s_2}-1}\Bigg)\nonumber\vsal
&=\bigg({\frac{\omega_q(s_2)^q}{s_2}-1}\bigg)-\frac{1}{q-1}\,\frac{p(p-1)}{p}\Bigg({\frac{\big({\tfrac{p}{p-1}}\big)^q}{s_2}-1}\Bigg)\nonumber\vsal
&=\bigg({\frac{\omega_q(s_2)^q}{s_2}-1}\bigg)-\frac{p-1}{q-1}\Bigg({\frac{\big({\tfrac{p}{p-1}}\big)^q}{s_2}-1}\Bigg)\,.\label{eq:3p9}
\end{align}

Since $s_2\geq H_q\big({\tfrac{p}{p-1}}\big)$ we have that $\omega_q(s_2)\leq \tfrac{p}{p-1}$, thus by \eqref{eq:3p9} $\gamma(s_1,s_2)\lessapprox 1-\frac{p-1}{q-1}=-\frac{p-q}{q-1}<0$. 

We prove now that for each $(s_1,s_2)$ which satisfies \eqref{eq:3p1} we have $\frac{\partial\, t}{\partial s_1}\leq 0$. Suppose that this result is not true. Then since for every $s_2$: $\gamma(s_1,s_2)\lessapprox 0$, as $s_1\to 0^{+}$, and $\delta(s_1,s_2)>0$, \;$\forall\,(s_1,s_2)$, we must have (by continuity reasons) that there exist $s_2\in(0,1)$,  $b\in(0,1)$, and  $\varepsilon>0$ such that: $\forall\,s_1\in[{b,b+\varepsilon})$ the quantity $t(s_1,s_2)$ is well defined, $(s_1,s_2)$ belongs to our domain of definition (see \eqref{eq:3p1}), and 
\begin{align}
\frac{\partial\,t}{\partial s_1}\cong \gamma(s_1,s_2)>0\,,\qquad\forall\,s_1\in({b,b+\varepsilon})\,,\label{eq:3p10}
\end{align}
while $\gamma(b,s_2)=0$. By \eqref{eq:3p10} and \eqref{eq:3p4} we get 
\begin{align}
\alpha(s_2)>\frac{1}{p(q-1)}\,\frac{(p-1)q\,\omega_q(\tau)-p(q-1)}{\omega_q(\tau)-1}\bigg({\frac{t^q}{s_2}-1}\bigg)\,,\;\forall\,s_1\in({b,b+\varepsilon})\,.\label{eq:3p11}
\end{align}
Also since  $\tau(s_1,s_2,t)=\tau=\frac{p-q}{p}\,\frac{t^p-s_1}{t^{p-q}-\frac{s_1}{s_2}}$, we have\footnote{We consider $t$ as a separate variable not depending on $s_1$, $s_2$.}

\begin{align*}
\frac{\partial\,\tau}{\partial\,t}&\cong p\,t^{p-1}\bigg({t^{p-q}-\frac{s_1}{s_2}}\bigg)-({t^p-s_1})(p-q)\,t^{p-q-1}\vsal
&\cong p\,t^{q}\bigg({t^{p-q}-\frac{s_1}{s_2}}\bigg)-(p-q)({t^p-s_1})\vsal
&=p\,t^{p}-p\,t^{q}\frac{s_1}{s_2}-(p-q)\,t^p+(p-q)s_1\vsal
&=q\,t^p-p\,t^q\frac{s_1}{s_2}+(p-q)s_1>0\,,
\end{align*}
as we had already proved.  Moreover if we consider $t$ as a fixed variable we have 
\begin{align*}
\frac{\partial\,\tau}{\partial s_1}=\frac{\partial}{\partial s_1}\tau(s_1,s_2,t)&\cong ({-1})\bigg({t^{p-q}-\frac{s_1}{s_2}}\bigg)-({t^p-s_1})\Big({-\frac{1}{s_2}}\Big)\vsal
  &=-t^{p-q}+\frac{t^p}{s_2}>0\,.
\end{align*}
Also since \eqref{eq:3p10} is true we have that: $\frac{\partial}{\partial s_1}t(s_1,s_2)>0$,\;$\forall\,s_1\in({b,b+\varepsilon})$. Thus if $s_1\in({b,b+\varepsilon})$ increases then $t(s_1,s_2)$ increases, and thus $\tau(s_1,s_2,t)$ increases. But \eqref{eq:3p11} can be written as 
\begin{align}
\alpha(s_2)>\frac{1}{p(q-1)}\bigg({(p-1)q+\frac{p-q}{\omega_q(\tau)-1}}\bigg)\bigg({\frac{t^q}{s_2}-1}\bigg)\,.\label{eq:3p12}
\end{align}
So, if we choose $s_1\in({b,b+\varepsilon})$ we must have $t(s_1,s_2):=t>t({b,s_2})$ and $\tau:=\tau(s_1,s_2)>\tau(b,s_2)\;\Rightarrow\; \frac{1}{\omega_q(\tau)-1}>\frac{1}{\omega_q(\tau(b,s_2))-1}$, thus \eqref{eq:3p12} gives 
\begin{align}
	\alpha(s_2)>\frac{1}{p(q-1)}\bigg({(p-1)q+\frac{p-q}{\omega_q\big(\tau(b,s_2)\big)-1}}\bigg)\bigg({\frac{t^q(b,s_2)}{s_2}-1}\bigg)\,,\label{eq:3p13}
\end{align}
where $\tau(b,s_2)=\frac{p-q}{p}\,\frac{t^p(b,s_2)-s_1}{t^{p-q}(b,s_2)-\frac{s_1}{s_2}}$\,. But the right hand side of \eqref{eq:3p13} equals $\alpha(s_2)$ since $\gamma(b,s_2)=0$. This gives a contradiction, thus $\frac{\partial\,t}{\partial s_1}\leq 0$, on the range of $(s_1,s_2)$, that \eqref{eq:3p1} describes.  Moreover we prove that there does not exist an interval $({\gamma,\gamma+\varepsilon})$ and an $s_2\in(0,1)$ such that $\frac{\partial}{\partial s_1}t(s_1,s_2)=0$\,,\;$\forall\,s_1\in ({\gamma,\gamma+\varepsilon})$. Indeed if such an interval exists, we should have $t(s_1,s_2)=c$\,,\;$\forall\,s_1\in ({\gamma,\gamma+\varepsilon})$, for some constant $c\geq 1$\,, and by \eqref{eq:3p4} 
\begin{align*}
\gamma(s_1,s_2)&=0\,\quad\forall\,s_1\in ({\gamma,\gamma+\varepsilon})\qquad\Rightarrow\vsal
	\alpha(s_2)& = \frac{1}{p(q-1)}\bigg({(p-1)q+\frac{p-q}{\omega_q({\tau})-1}}\bigg)\bigg({\frac{c^q}{s_2}-1}\bigg)
\end{align*}
that is $\tau$ depends only on $s_2$, which is a contradiction since for every $s_1\in ({\gamma,\gamma+\varepsilon})$ we then have $\tau(s_1,s_2,t)=\tau=\frac{p-q}{p}\,\frac{c^p-s_1}{c^{p-q}-\frac{s_1}{s_2}}$, which is obviously depends on $s_1$.  

\emph{We have proved the following theorem}

\begin{theorem}
The function $t(s_1,s_2)$ which is constructed in \cite{12} and studied in \cite{13}, is strictly decreasing on the variable $s_1$ (with $s_2$ being fixed).
\end{theorem}
\qed

\vspace{50pt}
\noindent Nikolidakis Eleftherios\\
Assistant Professor\\
Department of Mathematics \\
Panepistimioupolis, University of Ioannina, 45110\\
Greece\\
E-mail address: enikolid@uoi.gr


\newpage



\begin{thebibliography}{99}
	
\bibitem {1}
D. L. Burkholder.
\emph{Martingales and Fourier analysis in Banach spaces,}
	Probability and analysis (Varenna 1985), 61--108, Lecture Notes in Math., 1206, Springer, Berlin, 1986.
	
\bibitem {2}
D. L. Burkholder.
\emph{Explorations in martingale theory and its applications,}
	{\'E}cole d'{\'E}t{\'e} de Probabilit{\'e}s de Saint-Flour XIX—1989, Springer, Berlin, Heidelberg, 1991. 1-66.
	
	 	
\bibitem{3}
Delis, Anastasios D.; Nikolidakis, Eleftherios N.,
\emph{Sharp integral inequalities for the dyadic maximal operator and applications.},
	Math. Z. 291 (2019), No. 3-4,
	1197--1209.
	
	
	\bibitem {4}A. D. Melas.
	\emph{Sharp general local estimates for dyadic-like maximal operators and related Bellman functions,}
	Adv. in Math. 220 (2009), no. 2, 367--426

	
	\bibitem{5}
	A. D. Melas,
	\emph{The Bellman functions of dyadic-like maximal operators and related inequalities},
	Adv. in Math. 192 (2005), no. 2,
	310--340.
	
	
	\bibitem{6}
	E. N. Nikolidakis,
	\emph{Extremal problems related to maximal dyadic like operators},
	J. Math. Anal. Appl. 369 (2010), no. 1,
	377--385.
	
	\bibitem{7}
	E. N. Nikolidakis,
	\emph{Extremal sequences for the Bellman function of the dyadic maximal operator},
	Rev. Mat. Iberoam. 33 (2017), no. 2, 489--508.
	
	\bibitem{8}
	E. N. Nikolidakis,
	\emph{Extremal Sequences for the Bellman Function of the Dyadic Maximal Operator and Applications to the Hardy Operator},
	Canadian Journal of Mathematics 69, No.6, (2017), 1364-1384.
	
	\bibitem{9}
	E. N. Nikolidakis,
	\emph{Extremal sequences for the Bellman function of three variables of the dyadic maximal operator in relation to Kolmogorov's inequality},
	Trans. Am. Math. Soc. 372, No.9, (2019), 6315-6342.
	
	\bibitem{10}
	E. N. Nikolidakis,
	\emph{The geometry of the dyadic maximal operator},
	Rev. Mat. Iberoam. 30 (2014), no.4, 1397--1411.
	
	\bibitem{11}
	E. N. Nikolidakis,
	\emph{A multiparameter integral inequality for the dyadic maximal operator and applications},
	arXiv:1905.08091 (submitted)
	
	\bibitem{12}
	E. N. Nikolidakis,
	\emph{Sharp upper bounds for integral quantities related to the Bellman function of three variables of the dyadic maximal operator},
	arXiv:2312.05498
	
	\bibitem{13}
	E. N. Nikolidakis,
	\emph{On a sharp upper bound related to the Bellman function of three integral variables of the dyadic maximal operator},
	arXiv:2401.14821
	
	\bibitem{14}
	E. N. Nikolidakis, A. D. Melas,
	\emph{A sharp integral rearrangement inequality for the dyadic maximal operator and applications},
	Appl. Comput. Harmon. Anal., 38 (2015), no. 2, 242--261.
	
	\bibitem {15}
	E. N. Nikolidakis, A. D. Melas.
	\emph{Dyadic weights on $\mathbb{R}^{n}$ and reverse H\"{o}lder inequalities,}
	Stud. Math. 234 (2016), no.3, 281--290.
	
	\bibitem{16}
	L. Slavin, A. Stokolos, V. Vasyunin,
	\emph{Monge-Amp\`{e}re equations and Bellman functions: The dyadic maximal operator,}
	C. R. Math. Acad. Sci. Paris S\'{e}r. I. 346 (2008), no. 9-10,
	585--588.
	
	\bibitem{17}
	L. Slavin, A. Volberg,
	\emph{The $s$-function and the exponential integral},
	Topivs in harmonic analysis and ergodic theory, 215--228, Contemp. Math. 444.
	Amer. Math. Soc., Providence, RI, 2007.
	
	\bibitem{18}
	V. Vasyunin,
	\emph{The exact constant in the inverse H\"{o}lder inequality for Muckenhoupt weights},
	Algebra $i$ Analiz, 15 (2003),
	no. 1, 73--117; translation in St. Petersburg Math. J. 15 (2004), no. 1, 49--79
	
	\bibitem{19}
	V. Vasyunin, A. Volberg,
	\emph{The Bellman function for the simplest two weight inequality: an investigation of a particular case},
	Algebra i Analiz 18 (2006), no. 2, 24--56, translation in
	St. Petersburg Math. J., 18 (2007), no. 2 pp 201--222.
	
	\bibitem{20}
	V. Vasyunin, A. Volberg,
	\emph{Monge-Amp\`{e}re equation and Bellman optimization of Carleson embedding theorems, }Linear and complex analysis, 195--238,
	Amer. Math. Soc. Transl. Ser. 2, 226,
	Amer. Math. Sci., 63,  Providence, RI, 2009.
	
	\bibitem {21}
	G. Wang,
	\emph{ Sharp maximal inequalities for conditionally symmetric martingales and Brownian motion,}
	Proceedings of the American Mathematical Society 112 (1991): 579-586.
	
	
\end{thebibliography}
\end{document}